\documentclass[a4paper]{amsart}
\usepackage{amsmath,amssymb,amsfonts}
\usepackage{colordvi}
\usepackage[all]{xy}

\def\mpoint{\;.}
\def\mvirg{\;,}

\def\mpn{\medskip\par\noindent}

\def\mmpn{\vskip 1em minus 1em\par\noindent}

\def\sp{\bigskip\par}

\def\smp{\smallskip\par}

\def\CB{{\mathcal B}}
\def\CC{{\mathcal C}}

\def\CE{{\mathcal E}}

\def\CR{{\mathcal R}}

\def\Vect{\operatorname{Vect}\nolimits}

\def\Aut{\operatorname{Aut}\nolimits}

\def\op{^{op}}
\def\ls#1#2{{\,^{#1}\!#2}}

\def\pf{\par\bigskip\noindent{\bf Proof~: }}
\def\endpf{~\hfill\rlap{\hspace{-1ex}\raisebox{.5ex}{\framebox[1ex]{}}\sp}\bigskip\pagebreak[3]}

\makeatletter
\renewenvironment{enumerate}{\ifnum \@enumdepth >3 \@toodeep\else
       \advance\@enumdepth \@ne
       \edef\@enumctr{enum\romannumeral\the\@enumdepth}\list
       {\csname  label\@enumctr\endcsname}{\setlength{\topsep}{1ex}
\setlength{\itemsep}{0 pt}\usecounter
         {\@enumctr}\def\makelabel##1{\hss\llap{##1}}}\fi}{\endlist}

\renewenvironment{itemize}{\ifnum \@itemdepth >3 \@toodeep\else
\advance\@itemdepth \@ne
\edef\@itemitem{labelitem\romannumeral\the\@itemdepth}
\list{\csname\@itemitem\endcsname}{\setlength{\topsep}{1ex}\setlength
{\itemsep}{0pt}\def\makelabel##1{\hss\llap{##1}}}\fi}
{\endlist}
\def\@seccntformat#1{\csname the#1\endcsname.\quad}

\def\section{\pagebreak[3]\setcounter{prop}{0}\setcounter{equation}{0}\@startsection{section}{1}{\z@}{4ex plus  6ex}{2ex}{\center\reset@font \large\bf}}
\newcommand{\subsect}[1]{\medskip\par\noindent\pagebreak[3]\refstepcounter{subsection}\refstepcounter{prop}{\bf \thesection.\arabic{prop}.\ #1.\ }}
\makeatother

\def\theprop{\thesection.\arabic{prop}}

\newenvironment{enonce}[1]{\pagebreak[3]\refstepcounter{prop}\mmpn
{{\bf  \thesection.\arabic{prop}.\ #1.}}\begin{it} }{\end{it}\smp}
\def\thesection{\arabic{section}}

\newcommand{\result}[1]{\begin{enonce}{#1}}
\newcommand{\fresult}{\end{enonce}}

\newcommand{\mbigvee}[1]{\mathop{\bigvee}_{#1}\limits}

\begin{document}

\centerline{\Large\bf The monoid algebra of all relations on a finite set\footnote{This paper will never appear in this form. It is replaced by the following expanded article:\\
- The algebra of Boolean matrices, correspondence functors, and simplicity, submitted preprint, 2018.}}
\vspace{.5cm}
\centerline{Serge Bouc and Jacques Th\'evenaz}
\vspace{1cm}

\begin{footnotesize}
{\bf Abstract~:} 
We classify all the simple modules for the algebra of relations on a finite set, give their dimension, and find the dimension of the Jacobson radical of the algebra.

\bigskip\par

{\bf AMS Subject Classification~:} 06A06, 13C05, 16B50, 16D60, 16N20, 18B05, 18B10, 18B35, 18E05.\par
{\bf Key words~:} finite set, relation, simple module, Jacobson radical.
\end{footnotesize}

%%%%%%%% Introduction

\section{Introduction}
\noindent
J.A.~Green was a pioneer in the representation theory of algebras, in particular group algebras, but he also produced important work in the theory of monoids. We wish to contribute to this memorial issue by presenting some results on the $k$-algebra $\CR_X$ of a specific monoid, namely the monoid of all relations on a finite set~$X$, where $k$ is a field. This is an algebra of dimension $2^{n^2}$, where $n=|X|$, hence growing very fast in terms of~$n$.

In a previous work~\cite{BT1}, we described all simple modules for the algebra $\CE_X$ of essential relations on~$X$, which is a quotient of~$\CR_X$, but it was then not clear how to extend this result. We do so here and find all the simple modules for~$\CR_X$. They are classified by isomorphism classes of triples~$(E,R,V)$, where $E$ is finite set with $|E|\leq|X|$, $R$ is a partial order relation on~$E$, and $V$ is a simple module for the group algebra $k\Aut(E,R)$. When $E=X$, we recover the simple modules for the essential algebra~$\CE_X$.

We also give a formula for the dimension of each simple module and describe the action of every relation. From this, we deduce a formula for the dimension of the Jacobson radical of~$\CR_X$ (in characteristic zero). The formulas behave exponentially with respect to~$n$.

Apart from~\cite{BT1}, the other main ingredient for this work is our paper~\cite{BT2} about correspondence functors, in which all simple correspondence functors are classified and a formula is obtained for the dimension of each of their evaluations. It turns out that the evaluation at~$X$ of a simple functor is either zero or a simple module for the algebra~$\CR_X$. Conversely every simple module for the algebra~$\CR_X$ occurs as the evaluation at~$X$ of a simple functor. This provides a way to handle simple modules for the algebra~$\CR_X$ by using simple correspondence functors. This method sheds some new light on the structure of the algebra~$\CR_X$ and allows us to prove our results.

Throughout this paper, $X$ and $E$ denote finite sets and $k$ denotes a field.

%%%%%%%% Section 2

\section{Simple modules for the essential algebra} \label{Section-essential}
\bigskip
\noindent
A {\em correspondence} from a finite set $X$ to a finite set~$Y$ is a subset of $Y\times X$. The set of all correspondences from~$X$ to~$Y$ is written $\CC(Y,X)$. We reverse the order of $X$ and $Y$ in order to have later a left action of correspondences.
Correspondences can be composed as follows: if $R\subseteq Z\times Y$ and $S\subseteq Y\times X$, then
$$RS=\{ (z,x)\in Z\times X \mid \exists y\in Y, (z,y)\in R, (y,x)\in S \} \subseteq Z\times X \mpoint$$
A correspondence from a finite set $X$ to itself is called a {\em relation} on~$X$.
The set $\CC(X,X)$ of all relations on~$X$ is a monoid for the composition above and forms a basis of the finite-dimensional $k$-algebra $\CR_X:=k\CC(X,X)$. This algebra is called the algebra of all relations on~$X$.

A relation~$R$ on~$X$ is called {\em essential\/} if it does not factor through a set of cardinality strictly smaller than~$|X|$.
The $k$-submodule generated by the set of inessential relations is a two-sided ideal
$$I_X=\sum_{|Y|<|X|} k\CC(X,Y)k\CC(Y,X)$$
and the quotient
$$\CE_X:=k\CC(X,X)/I_X$$
is called the {\em essential algebra}.

By an order relation, we always mean a partial order. Let $R$ be an order relation on a finite set~$E$ and let $\sigma\in\Sigma_E$, the group of all permutations of the set~$E$. Define the conjugate relation $\ls \sigma R=\Delta_\sigma R \Delta_{\sigma^{-1}}$, where
$$\Delta_\sigma=\{(\sigma(x),x) \in E\times E\mid x\in E\} \mpoint$$
Then $\ls \sigma R$ is an order relation on~$E$ and this defines an action of the group $\Sigma_E$ on the set of all order relations on~$E$. The stabilizer of $R$ under this action is written $\Aut(E,R)$. 
The group $\Sigma_E$ acts by conjugation on the set of all pairs $(R,V)$, where $R$ is an order relation on~$E$ and $V$ is a simple $k\Aut(E,R)$-module. The following result is Theorem~8.1 in~\cite{BT1}.

\result{Theorem} \label{simple-CE} Let $E$ be a finite set.
The set of isomorphism classes of simple $\CE_E$-modules is parametrized by the set of conjugacy classes of pairs $(R,V)$, where $R$ is an order relation on~$E$ and $V$ is a simple $k\Aut(E,R)$-module.
\fresult

We let $T_{R,V}$ be the simple $\CE_E$-module parametrized by (the conjugacy class of) the pair $(R,V)$.

%%%%%%%% Section 3

\section{Simple correspondence functors} \label{Section-functors}

\bigskip
\noindent
We define the category~$\CC$ as follows~:
\begin{itemize}
\item The objects of $\CC$ are the finite sets.
\item For any two objects $X$ and $Y$, the set of morphisms from $X$ to~$Y$ is the set $\CC(Y,X)$ of all correspondences from $X$ to~$Y$, namely all subsets of $Y\times X$.
\item The composition of morphisms is defined as above.
\end{itemize}
The {\em $k$-linearization} of the category~$\CC$ is defined as follows~:
\begin{itemize}
\item The objects of $k\CC$ are the objects of $\CC$.
\item For any two objects $X$ and $Y$, the set of morphisms from $X$ to~$Y$ is the
$k$-vector space $k\CC(Y,X)$ with basis $\CC(Y,X)$.
\item The composition of morphisms in~$k\CC$ is the $k$-bilinear extension
$$k\CC(Z,Y) \times k\CC(Y,X) \longrightarrow k\CC(Z,X)$$
of the composition in~$\CC$.
\end{itemize}

A {\em correspondence functor} is a $k$-linear functor from $k\CC$ to the category $k\text{-\!}\Vect$ of $k$-vector spaces.
We could define a correspondence functor to be a functor from $\CC$ to $k\text{-\!}\Vect$, but it is convenient to linearize first the category~$\CC$ (just as for group representations, where one can first introduce the group algebra).
The category of correspondence functors is an abelian category and we are particularly interested in its simple objects.

We now describe the parametrization of simple correspondence functors. If $S$ is a simple correspondence functor, we let $E$ be a minimal set such that $S(E)\neq\{0\}$. Then $S(E)$ is a module for the algebra $\CR_E=k\CC(E,E)$, and actually a module for the essential algebra $\CE_E$, because the ideal $I_E$ of inessential relations acts by zero, by minimality of~$E$. Moreover, $S(E)$ is easily seen to be simple (see Proposition~2.6 in~\cite{BT2}). Therefore, by Theorem~\ref{simple-CE}, $S(E)\cong T_{R,V}$ for some order relation~$R$ on~$E$ and some simple $k\Aut(E,R)$-module~$V$. Associated to the simple functor $S$, this defines a triple $(E,R,V)$ which is uniquely defined up to isomorphism. The following result is Theorem~3.12 in~\cite{BT2}.

\result{Theorem} \label{simple-functors} 
The set of isomorphism classes of simple correspondence functors is parametrized by the set of isomorphism classes of triples $(E,R,V)$, where $E$ is a finite set, $R$ is an order relation on~$E$, and $V$ is a simple $k\Aut(E,R)$-module.
\fresult

We let $S_{E,R,V}$ be the simple correspondence functor parametrized by (the isomorphism class of) the triple $(E,R,V)$. Thus we have
$$S_{E,R,V}(E)\cong T_{R,V} \quad \text{ and } \quad S_{E,R,V}(F)=\{0\} \;\text{ if } |F|<|E| \mpoint$$
One of the main results of~\cite{BT2} provides a closed formula for the dimension of any evaluation of a simple functor.
The following result is Theorem~17.22 in~\cite{BT2}.

\result{Theorem} \label{dim-simple} Let $S_{E,R,V}$ be the simple correspondence functor parametrized by the triple $(E,R,V)$ and let $X$ be a finite set. Then
$$\dim S_{E,R,V}(X)=\frac{\dim V}{|\Aut(E,R)|} \sum_{i=0}^{|E|}(-1)^i\binom{|E|}{i}(g_{E,R}-i)^{|X|} \mvirg$$
where $g_{E,R}$ is a positive integer canonically associated to~$(E,R)$.
\fresult

We only sketch here where the integer $g_{E,R}$ comes from and give references. We can associate to the poset $(E,R)$ a finite lattice $T$ having $(E,R\op)$ as the full subset of its join-irreducible elements, where $R\op$ denotes the opposite relation. For instance, we can take $T$ to be the lattice of all lower ideals in the poset $(E,R\op)$ and identify $E$ with the set of principal ideals. Then there are two operations $s^\infty$ and $r^\infty$ in~$T$, defined in Section~16 of~\cite{BT2}, and we let
$$G_{E,R}=E\sqcup \{a\in T\mid r^\infty(s^\infty(a))=a \} \mpoint$$
Clearly $E\subseteq G_{E,R} \subseteq T$ and we let $g_{E,R}=|G_{E,R}|$. Since $\dim S_{E,R,V}(X)$ does not depend on the choice of~$T$, the integer $g_{E,R}$ turns out to be independent of the choice of~$T$. We refer to Section~16 of~\cite{BT2} for details.

%%%%%%%% Section 4

\section{Simple modules for the algebra of relations} \label{Section-algebra}

\bigskip
\noindent
We now give the parametrization of all simple modules for the algebra $\CR_X$ of relations on a finite set~$X$. They are described in terms of simple correspondence functors.

\result{Theorem} \label{simple-RX} Let $X$ be a finite set.
\begin{enumerate}
\item The set of isomorphism classes of simple $\CR_X$-modules is parametrized by the set of isomorphism classes of triples $(E,R,V)$, where $E$ is a finite set with $|E|\leq |X|$, $R$ is an order relation on~$E$, and $V$ is a simple $k\Aut(E,R)$-module.
\item The simple module parametrized by the triple $(E,R,V)$ is $S_{E,R,V}(X)$, where $S_{E,R,V}$ is the simple correspondence functor corresponding to the triple $(E,R,V)$.
\end{enumerate}
\fresult

\pf
We first claim that any simple $\CR_X$-module~$S$ occurs as the evaluation of some simple functor $S_{E,R,V}$, that is, $S\cong S_{E,R,V}(X)$. This forces $|E|\leq|X|$ because $S_{E,R,V}$ vanishes on sets~$Y$ with $|E|>|Y|$. A proof of this first claim appears in Proposition~3.2 of~\cite{We} and is attributed to Green (6.2 in~\cite{Gr}). This requires to view $\CR_X=k\CC(X,X)$ as a category with a single object~$X$, hence a full subcategory of~$k\CC$. Another proof of the claim appears in Corollary~4.2.4 of~\cite{Bo}, but only for subcategories of the biset category. However, the result and its proof hold for any small category.\par

Our second claim is that, conversely, the evaluation of a simple functor $S_{E,R,V}$ at a finite set~$X$ is either zero or a simple $\CR_X$-module. A proof of this second claim again appears in Proposition~3.2 of~\cite{We}, or in Corollary~4.2.4 of~\cite{Bo} (for subcategories of the biset category, but the result is again quite general). Moreover, $S_{E,R,V}(X)$ is zero if $|E|>|X|$ and is nonzero if $|E|\leq|X|$ because $S_{E,R,V}(E)=T_{R,V}$ is nonzero, hence $S_{E,R,V}(X)$ is also nonzero, by Corollary~4.4 in~\cite{BT2}.\par

The given references also indicate that the simple functor $S_{E,R,V}$ associated to a simple $\CR_X$-module~$S$, i.e. such that $S\cong S_{E,R,V}(X)$, is unique up to isomorphism. This provides the parametrization of the statement and completes the proof.
\endpf

Note that we used in the proof above the non-vanishing of $S_{E,R,V}(X)$ when $|E|\leq|X|$. This is a special property of correspondence functors (Corollary~4.4 in~\cite{BT2}) and it may not hold for representations of other small categories. For instance, it does not hold for the biset category (see~\cite{BST}).\par

In view of Theorem~\ref{simple-RX}, a formula for the dimension of any simple $\CR_X$-module is now given by Theorem~\ref{dim-simple}.\par

An explicit description can be given for every simple $\CR_X$-module and for the action of relations. This is rather technical and is an application of Theorem~18.4 in~\cite{BT2}. It can be summarized as follows. We fix a simple $\CR_X$-module $S_{E,R,V}(X)$. As already noticed, we can associate to the poset $(E,R)$ a finite lattice $T$ having $(E,R\op)$ as the full subset of its join-irreducible elements, and we can also choose $T$ such that $\Aut(T)=\Aut(E,R\op)$ by taking for instance $T$ to be the lattice of all lower ideals in the poset $(E,R\op)$. Then there is a correspondence functor $F_T$ defined by $F_T(X)=k(T^X)$, the $k$-vector space with basis the set $T^X$ of all maps from $X$ to~$T$. The action of a correspondence $U\in\CC(Y,X)$ on a map $\varphi\in T^X$ is a map $U\varphi \in T^Y$ defined by
$$(U\varphi)(y)=\mbigvee{(y,x)\in U}\varphi(x)\mvirg$$
and this makes $F_T$ into a correspondence functor, by Proposition~11.2 in~\cite{BT2}. Define a subset
$$\CB_{E,R,X}=\{\varphi \in T^X \mid E\subseteq \varphi(X) \subseteq G_{E,R} \} \mvirg$$
where $G_{E,R}$ is defined as in Section~\ref{Section-functors}. Note that the cardinality of $\CB_{E,R,X}$ is given by the formula
$$|\CB_{E,R,X}|=\sum_{i=0}^{|E|}(-1)^i\binom{|E|}{i}(g_{E,R}-i)^{|X|} \mpoint$$
This gives an interpretation of the sum appearing in Theorem~\ref{dim-simple}. The formula is proved in Lemma~8.1 of~\cite{BT2}.
Let $k\CB_{E,R,X}$ be the $k$-subspace of $F_T(X)$ with basis~$\CB_{E,R,X}$. The family of subspaces $k\CB_{E,R,X}$ do not form a subfunctor of~$F_T$, but they can be used to describe the evaluations of the simple functor $S_{E,R,V}$. Indeed, by Remark 4.8 and Corollary~17.17 in~\cite{BT2}, we have
$$S_{E,R,V}(X) \cong k\CB_{E,R,X} \otimes_{k\Aut(E,R)} V \mvirg$$
using the right action of $\Aut(E,R)$ on~$\CB_{E,R,X}$ defined as follows: given $\varphi\in\CB_{E,R,X}$ and a permutation $\sigma\in\Aut(E,R)$, then $\sigma$ yields an automorphism of~$T$ since $\Aut(E,R\op)=\Aut(E,R)$ and we let $\varphi\cdot\sigma:= \sigma^{-1}\circ \varphi$. This provides an explicit description of the simple $\CR_X$-module $S_{E,R,V}(X)$ as a $k$-vector space. It remains to describe the action of every relation $U\in\CC(X,X)$.\par

Let $\pi_{T,X}$ be the $k$-linear idempotent endomorphism of $k(T^X)$ defined by
$$\forall \varphi\in T^X,\;\;\pi_{T,X}(\varphi)=\left\{\begin{array}{ll}\varphi&\hbox{if}\;E \subseteq \varphi(X) \mvirg\\
0&\hbox{otherwise} \mpoint\end{array}\right.$$
By Theorem~17.7 in~\cite{BT2}, there is an element $u_T\in k(T^T)$ which has the property that, for any $\varphi\in T^X$, the composition $u_T\circ \varphi$ is a $k$-linear combination of maps $\psi\in T^X$ such that $\psi(X)\subseteq G_{E,R}$. (Actually, $u_T$ is idempotent, by Theorem~17.9 in~\cite{BT2}.)
Then, for any map $\varphi\in T^X$, we obtain
$$\pi_{T,X}(u_T\circ \varphi) \in k\CB_{E,R,X} \mvirg$$
that is, a $k$-linear combination of maps $\psi\in T^X$ such that $E\subseteq \psi(X)\subseteq G_{E,R}$.
Now we can describe the action of relations on the simple $\CR_X$-module $S_{E,R,V}(X)$.

\result{Theorem} \label{action} Fix the notation above.
\begin{enumerate}
\item $S_{E,R,V}(X)\cong k\CB_{E,R,X} \otimes_{k\Aut(E,R)} V$ as $k$-vector spaces.
\item The action of a relation $U\in\CC(X,X)$ on an element
$$\varphi\otimes v \in k\CB_{E,R,X} \otimes_{k\Aut(E,R)} V \,, \qquad(\varphi\in\CB_{E,R,X}, \; v\in V)$$
is given by
$$U\cdot (\varphi\otimes v)=\pi_{T,X}(u_T\circ U\varphi) \otimes v \mpoint$$
\end{enumerate}
\fresult

The proof appears in Theorem~18.4 of~\cite{BT2}. Theorem~\ref{action} provides a computational method for the description of any simple $\CR_X$-module.

%%%%%%%% Section 5

\section{The Jacobson radical of the algebra of relations} \label{Section-Jacobson}

\bigskip
\noindent
In this section, we assume for simplicity that the field $k$ has characteristic zero. Our purpose is to give the dimension of the Jacobson radical $J(\CR_X)$ of the $k$-algebra~$\CR_X$.

\result{Theorem} \label{Jacobson}
Assume that $k$ is a field of characteristic zero.
Let $J(\CR_X)$ be the Jacobson radical of the $k$-algebra~$\CR_X$ and let $n=|X|$. Then
$$\dim J(\CR_X) = 2^{n^2} - \sum_{e=0}^n \sum_R \displaystyle\frac{1}{|\Aut(E,R)|}
\Big(\sum_{i=0}^e (-1)^i \binom{e}{i}(g_{E,R}-i)^{n}  \Big)^2 \mvirg$$
where $R$ runs over a set of representatives of $\Sigma_e$-conjugacy classes of order relations on the set~$E=\{1,\ldots,e\}$. The integer $g_{E,R}$ is as in Theorem~\ref{dim-simple}.
\fresult

\pf
Since $k$ has characteristic zero, the semi-simple algebra $\CR_X/J(\CR_X)$ is separable, that is, it remains semi-simple after scalar extension to an algebraic closure $\overline k$ of~$k$. In other words, $\dim J(\CR_X)$ does not change after this scalar extension. Therefore, we can assume that $k=\overline k$.\par

By Theorem~\ref{simple-RX}, every simple $\CR_X$-module has the form $S_{E,R,V}(X)$ with $|E|\leq|X|$, where $S_{E,R,V}$ is the simple correspondence functor parametrized by the triple $(E,R,V)$. In order to have a parametrization, we take $E=\{1,\ldots,e\}$ with ${0\leq e\leq n}$, we take $R$ in a set of representatives as in the statement, and finally we take $V$ in a set of representatives of isomorphism classes of simple $k\Aut(E,R)$-modules.\par

Since the endomorphism algebra of a simple module is isomorphic to~$k$, by Schur's lemma and the assumption that $k$ is algebraically closed, the dimension of the semi-simple algebra $\CR_X/J(\CR_X)$ is equal to the sum of the squares of the dimensions of all simple modules, by Wedderburn's theorem. By Theorem~\ref{dim-simple}, we obtain
$$\begin{array}{rcl}
\dim \big(\CR_X/J(\CR_X)\big) &=& \displaystyle \sum_{E,R,V} \big(\dim S_{E,R,V}(X)\big)^2 \\
&=& \displaystyle \sum_{E,R,V} \Big(\frac{\dim V}{|\Aut(E,R)|}\Big)^2 \Big(\sum_{i=0}^{|E|}(-1)^i\binom{|E|}{i}(g_{E,R}-i)^{|X|}\Big)^2 \\
&=& \displaystyle \sum_{e=0}^n \sum_R \Big(\sum_V \frac{(\dim V)^2}{|\Aut(E,R)|^2}\Big)
\Big(\sum_{i=0}^e (-1)^i \binom{e}{i}(g_{E,R}-i)^{n}\Big)^2 \\
&=& \displaystyle \sum_{e=0}^n \sum_R \frac{1}{|\Aut(E,R)|}
\Big(\sum_{i=0}^e (-1)^i \binom{e}{i}(g_{E,R}-i)^{n}\Big)^2
\mvirg\end{array}$$
because $\displaystyle\sum_V(\dim V)^2=\dim (k\Aut(E,R)) =|\Aut(E,R)|$, by semi-simplicity of the group algebra in characteristic zero (Maschke's theorem). Now
$$\dim J(\CR_X) = \dim \CR_X- \dim \big(\CR_X/J(\CR_X)\big) = 2^{n^2} - \dim \big(\CR_X/J(\CR_X)\big)$$
and the result follows.
\endpf

If $k$ is an algebraically closed field of prime characteristic~$p$, the formula has to be modified in a straightforward manner, in order to take into account the Jacobson radical of $k\Aut(E,R)$.
Then it seems likely that the same formula holds over any field of characteristic~$p$ (that is, $\CR_X/J(\CR_X)$ is likely to be a separable algebra), but we leave this question open.

%%%%%%%% Section 6

\section{Examples} \label{Section-examples}

\bigskip
\noindent
In this final section, we provide examples of the algebra $\CR_X$ for small values of $n=|X|$.

\subsect{Example} \label{empty}
Let $X=\emptyset$, hence $n=0$. There is a single relation on~$\emptyset$, namely~$\emptyset$, and $\CR_\emptyset\cong k$. The relation $\emptyset$ has a trivial automorphism group, which has only the trivial module~$k$ as a simple module, and we obtain that $\CB_{\emptyset,\emptyset,\emptyset}$ has cardinality~1 and $k\CB_{\emptyset,\emptyset,\emptyset}\cong k$, so that
$$S_{\emptyset,\emptyset,k}(\emptyset)\cong k\CB_{\emptyset,\emptyset,\emptyset} \otimes_k k\cong k \mpoint$$
The unique relation $\emptyset$ acts as the identity on~$k$.

\bigskip

\subsect{Example} \label{one}
Let $X=\{1\}$, hence $n=1$. There are 2 relations on~$\{1\}$, namely $\emptyset$ and $R=\{1\}\times \{1\}$, so $\CR_{\{1\}}$ has dimension~2.\par

For $E=\emptyset$, there is a simple correspondence functor $S_{\emptyset,\emptyset,k}$ and again $\CB_{\emptyset,\emptyset,\{1\}}$ has cardinality~1. Thus we get a simple $\CR_{\{1\}}$-module
$$S_{\emptyset,\emptyset,k}(\{1\})\cong k\CB_{\emptyset,\emptyset,\{1\}} \otimes_k k\cong k \mpoint$$
Both relations act as the identity on~$k$.\par

For $E=\{1\}$, there is a single essential relation~$R$, with a trivial automorphism group, having only the trivial module~$k$ as a simple module. It turns out that $g_{E,R}=2$ and $|\CB_{\{1\},R,\{1\}}|=2^1-1^1=1$.
Then there is a simple correspondence functor $S_{\{1\},R,k}$ and we have a simple $\CR_{\{1\}}$-module
$$S_{\{1\},R,k}(\{1\})\cong k\CB_{\{1\},R,\{1\}} \otimes_k k \cong k \mpoint$$
The relation $\emptyset$ acts by zero, while $R$ acts as the identity.\par

Since $\dim \CR_{\{1\}}=2$ and there are two simple modules, we have $\CR_{\{1\}}\cong k\times k$, a semi-simple algebra.

\bigskip

\subsect{Example} \label{two}
Let $X=\{1,2\}$, hence $n=2$. There are $2^4=16$ relations on~$\{1,2\}$, so $\CR_{\{1,2\}}$ has dimension~16.\par

For $E=\emptyset$, we get a simple $\CR_{\{1,2\}}$-module $S_{\emptyset,\emptyset,k}(\{1,2\})$, of dimension~1.\par

For $E=\{1\}$, there is a single essential relation~$R$, with a trivial automorphism group, having only the trivial module~$k$ as a simple module. It turns out that $g_{\{1\},R}=2$ and $|\CB_{\{1\},R,\{1,2\}}|=2^2-1^2=3$. We get a simple $\CR_{\{1,2\}}$-module $S_{\{1\},R,k}(\{1,2\})$ of dimension~3.\par

For $E=\{1,2\}$, there are two essential relations up to conjugacy, namely the equality relation~$\rm{eq}$ and the usual total order $\rm{tot}$. It turns out that $g_{\{1,2\},\rm{eq}} = 4$ and $\CB_{\{1,2\},\rm{eq},\{1,2\}}$ has cardinality
$$|\CB_{\{1,2\},\rm{eq},\{1,2\}}|=\sum_{i=0}^2 (-1)^i \binom{2}{i}(g_{\{1,2\},\rm{eq}}-i)^2=4^2-2\cdot 3^2 +2^2=2 \mpoint$$
Moreover, $\Aut(\{1,2\},\rm{eq})=C_2$ has order~2, with two simple modules $k_+$ and $k_-$ (assuming that the characteristic of $k$ is not~2). Therefore, we obtain two simple $\CR_{\{1,2\}}$-modules
$$\begin{array}{rcl}
S_{\{1,2\},\rm{eq},k_+}(\{1,2\}) &\cong& k\CB_{\{1,2\},\rm{eq},\{1,2\}} \otimes_{kC_2} k_+ \mvirg \\
S_{\{1,2\},\rm{eq},k_-}(\{1,2\}) &\cong& k\CB_{\{1,2\},\rm{eq},\{1,2\}} \otimes_{kC_2} k_- \mvirg
\end{array}
$$
both of dimension~1.\par

For the other relation $\rm{tot}$, we have $g_{\{1,2\},\rm{tot}} = 3$ and $\CB_{\{1,2\},\rm{tot},\{1,2\}}$ has cardinality
$$|\CB_{\{1,2\},\rm{tot},\{1,2\}}|=\sum_{i=0}^2 (-1)^i \binom{2}{i}(g_{\{1,2\},\rm{tot}}-i)^2=3^2-2\cdot 2^2 +1^2=2 \mpoint$$
Moreover, $\Aut(\{1,2\},\rm{tot})$ is the trivial group and has only the trivial module~$k$ as a simple module.
We obtain a simple $\CR_{\{1,2\}}$-module
$$S_{\{1,2\},\rm{tot},k}(\{1,2\}) \cong k\CB_{\{1,2\},\rm{tot},\{1,2\}} \otimes_k k\cong k\CB_{\{1,2\},\rm{tot},\{1,2\}}
\mvirg$$
of dimension~2.\par

Altogether, there are 5 simple $\CR_{\{1,2\}}$-modules and the dimension of the semi-simple quotient is
$$\dim \big(\CR_{\{1,2\}}/J(\CR_{\{1,2\}})\big)=1^2+3^2+1^2+1^2+2^2=16 \mvirg$$
so the Jacobson radical has dimension $\dim J(\CR_{\{1,2\}})=2^4-16=0$. Therefore $\CR_{\{1,2\}}$ is semi-simple (provided the characteristic of $k$ is not~2).

\bigskip

\subsect{Example} \label{more}
\def\posetVop{^{\mathop{\mathop{\displaystyle\bullet\;\bullet}^{\displaystyle/\backslash}\limits}^{\raisebox{-1ex}{$\displaystyle\bullet$}}\limits}}
\def\posetV{\raisebox{-2ex}{$^{\,\mathop{\mathop{\displaystyle\bullet\;\bullet}_{\displaystyle\backslash/}\limits}_{^{\displaystyle\bullet}}\limits}$}}
\def\ipoint{\mathop{\rule{.1ex}{.9ex}\hspace{.2ex}}^{_{\displaystyle\bullet}}_{^{\displaystyle\bullet}}\limits{\displaystyle\bullet}}
\def\totdeux{\mathop{\rule{.1ex}{.9ex}\hspace{.2ex}}\limits^{_{\displaystyle\bullet}}_{^{\displaystyle\bullet}}}
\def\tottrois{\mathop{\totdeux}\limits^{\mathop{\rule{.1ex}{.9ex}}\limits^{\displaystyle\bullet}}}
For $n=|X|\geq 3$, the algebra $\CR_X$ is not anymore semi-simple. Using the computer software GAP (\cite{GAP4}), we have computed the dimension of the Jacobson radical of $\CR_X$ for a set $X$ of cardinality 3, and it is equal to~42. According to Theorem~\ref{Jacobson}, this value can be recovered as follows:
$$\begin{array}{|c|c|c|c|c|r|}
\hline
{\rm Size}\;e&{\rm Poset}(E,R)&|\Aut(E,R)|&g_{E,R}&\sum\limits_{i=0}^e(-1)^i\binom{e}{i}(g_{E,R}-i)^3&{\rm total}\\
\hline
0&\emptyset&1&1&1&1\\
\hline
1&\bullet&1&2&7&49\\
\hline
2&\bullet\bullet&2&4&18&162\\
&\rule{0ex}{4ex}\totdeux&1&3&12&144\\
\hline
3&\bullet\bullet\bullet&6&5&6&6\\
&\rule{0ex}{4ex}\ipoint&1&5&6&36\\
&\rule{0ex}{5ex}\posetV&2&5&6&18\\
&\rule{0ex}{5ex}\raisebox{-2ex}{$\posetVop$}&2&5&6&18\\
&\tottrois&1&6&36&36\\
\hline
\end{array}
$$
In this case, the algebra $\CR_X$ has dimension $2^{3^2}=512$. The sum of the last column of this table is equal to 470, so we recover the dimension of the radical 42=512-470.\vspace{1ex}\par

For $|X|\geq 4$, the algebra $\CR_X$ has dimension greater than or equal to $2^{4^2}=65,\!536$. The direct computation of the radical of such a big algebra seems out of reach of usual computers. However, using the formula of Theorem~\ref{Jacobson} and the structure of the 16 posets of cardinality 4, one can show that the radical of $\CR_X$ has dimension $32,616$ when $|X|=4$.

\bigskip
\bigskip
\noindent
Serge Bouc, CNRS-LAMFA, Universit\'e de Picardie - Jules Verne,\\
33, rue St Leu, F-80039 Amiens Cedex~1, France.\\
{\tt serge.bouc@u-picardie.fr}

\medskip
\noindent
Jacques Th\'evenaz, Section de math\'ematiques, EPFL, \\
Station~8, CH-1015 Lausanne, Switzerland.\\
{\tt Jacques.Thevenaz@epfl.ch}

\end{document}